\documentclass{amsart}
\usepackage{docmute}
\usepackage{graphicx}
\usepackage{latexsym}
\usepackage{color}
\usepackage{amscd}
\usepackage[all]{xy}
\usepackage{enumerate}
\usepackage{subfigure}
\usepackage{soul}
\usepackage{comment}
\newtheorem{dfn}{Definition}[section]
\newtheorem{thm}[dfn]{Theorem}

\newtheorem{lem}[dfn]{Lemma}
\newtheorem{rem}[dfn]{Remark}

\newtheorem{prob}[dfn]{Problem}

\def\CM{\mathcal M}
\def\IT{\mathcal T}

\newcommand{\R}{\mathbb{R}}

\newcommand{\Z}{\mathbb{Z}}

\begin{document}

%----タイトル----%
\title[Generating twist subgroup]{Generating twist subgroup of mapping class group of non-orientable surface by involutions}

%----著者----%
\author{Kazuya Yoshihara}
\date{}

%----アブストラクト----%
\begin{abstract}
    Let $N_{g}$ denote the closed non-orientable surface of genus $g$ and let $\CM _g$ denote the mapping class group of $N_{g}$. Let $\IT _g$ denote the twist subgroup of $\CM _g$ which is the subgroup of $\CM _g$ is generated by all Dehn twists.
    %Du found the generating set for $\IT _g$ consisting of three elements of finite order in the case of odd genus. It is not known the generating set for $\IT _g$ consisting only of involutions.
    In this thesis, we proved that $\IT _g$ is generated by six involutions for $g \geq 16$ or $g=14$.
\end{abstract}
\maketitle

%----§1.Introduction----%
\section{Introduction}
\label{section:one}
\begin{figure}[ht]
    \begin{center}
        \includegraphics[scale=1.0,clip]{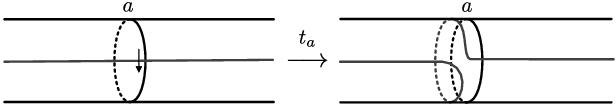}
    \end{center}
    \caption{Dehn twist along a simple closed curve $a$}
    \label{fig:dehn_twist}
\end{figure}

%点付き有向閉曲面の写像類群の導入%
For $g \geq 1$, let $N_{g}$ denote the closed connected non-orientable surface of genus $g$. The {\it{mapping class group}} $\CM _g$ is the group of isotopy classes of diffeomorphisms of $N_{g}$.
A simple closed curve $c$ in $N_{g}$ is {\it two-sided} (resp.\ one-sided) if a regular neighborhood of $c$, denoted by $N_c$, is an annulus (resp.\ a M\"obius band). Let $a$ be a two-sided simple closed curve on $N_{g}$. By the definition, the regular neighborhood of $a$ is an annulus, and it has two possible orientation. Now, we fix one of its two possible orientations. For two sided simple closed curve $a$, we can define the Dehn twist $t_a$. We indicate the direction of a Dehn twist by an arrow beside the curve $a$ as shown in Fig.~\ref{fig:dehn_twist}.

%TODO: Twist Subgroupの歴史が触れられていない 
For $g\geq2$, Lickorish (\cite{li1}, \cite{li2}) first proved that $ \CM _g$ is not generated by Dehn twists, and $ \CM _g$ is generated by Dehn twists and Y-homeomorphisms. Chillingworth \cite{chil} found a finite set of generators of this group. The number of Chillingworth's generators is improved to $g+1$ by Szepietowski \cite{sez3}. Hirose \cite{hirose} proved that his generating set is the minimal generating set by Dehn twists and Y-homemorphisms.
%Szepietowski \cite{sez1} proved that $ \CM _g$ is generated by involutions. The cardinality of his generating set of involutions depends linearly on $g$. 
Szepietowski show that $ \CM _g$ is generated by four involutions for $g \geq 4$ \cite{sez2}.
Denote by $\IT _g$ the subgroup of $ \CM _g$ is generated by all Dehn twists. We call $\IT _g$ the {\it{twist subgroup}} of $\CM _g$.
The group $\IT _g$ is an index $2$ subgroup of $ \CM _g$ (see \cite{li1}). In particular, $\IT _g$ is finitely generated.
Chillingworth \cite{chil} proved that $\IT _g$ is generated by a Dehn twist if $g=2$, two Dehn twists if $g=3$, $(3g-1)/2$ Dehn twists if $g\geq5$ is odd, $(3g)/2$ Dehn twists if $g\geq4$ is even. Stukow \cite{st2} proved that $\IT _g$ admits a finite presentation. Omori \cite{o} showed that $\IT _g$ is generated by $g+1$ Dehn twists for $g\geq4$. Recently, Du \cite{du} showed that $\IT _g$ is generated by three elements of finite order if $g \geq 5$ is odd. One of his generator is of order $2g$ and the other two are of order $2$.
%Duの結果を書く
%\begin{thm}
%\label{thm:du}
%If $g \geq 5$ is odd, $\IT _g$ is generated by three elements of finite order. One of the generator is of order $2g$. The other two are of order $2$.
%\end{thm}

On the other hand, it is not known the generating set for $\IT _g$ consisting only of involutions. In this paper, we showed that $\IT _g$ is generated by involutions.
\begin{thm}
    \label{thm:main}
    For $g \geq 16$ or $g=14$, $\IT _g$ is generated by six involutions.
\end{thm}

The paper is organized as follows.
In Section $2$ we recall the properties of Dehn twists and Y-homeomorphisms. In Section $3$ we construct involutions of $ \IT _g$ and prove the theorem~\ref{thm:main}. %Finally, in Section $4$, .
%----§1. Introduction ここまで----%

%----§2. Preliminaries ----%
\section{Preliminaries}
\label{section:two}
\begin{figure}[ht]
    \begin{center}
        \includegraphics[scale=1.0,clip]{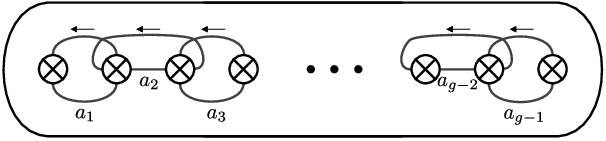}
    \end{center}
    \caption{Surface $N_{g}$ and its simple closed curves}
    \label{fig:generator1}
\end{figure}

%曲線の定義
We represent the surface $N_{g}$ as a connected sum of $g$ projective planes as in Fig.~\ref{fig:generator1}. In this Figure, each encircled cross mark represents a crosscap: the interior of the encircled disk is to be removed and each pair of antipodal points on the boundary are to be identified.

\begin{figure}[ht]
    \begin{center}
        \includegraphics[scale=1.0,clip]{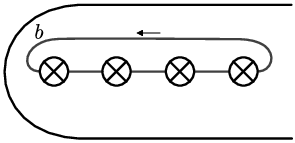}
        \hspace{0.4cm}
        \includegraphics[scale=1.0,clip]{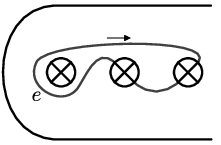}
    \end{center}
    \caption{Simple closed curves $b$ and $e$}
    \label{fig:curve_b_e}
\end{figure}

Let $a_1,a_2,\ldots,a_{g-1},b,$ and $e$ be two-sided simple closed curves on $N_g$ as in Fig.~\ref{fig:generator1} and Fig.~\ref{fig:curve_b_e}.

%\begin{figure}
%\begin{center}
%\includegraphics[scale=1.0,clip]{./drawable/generators3}
%\end{center}
%\caption{The simple closed curve $e$}
%\label{fig:curve_e}
%\end{figure}

\begin{figure}[ht]
    \begin{center}
        \includegraphics[scale=1.0,clip]{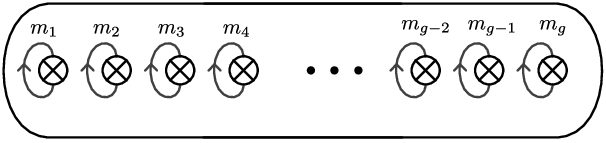}
    \end{center}
    \caption{One-sided simple closed curves $m_1,m_2,\ldots,m_{g}$}
    \label{fig:one_sided_simple_closed_curves}
\end{figure}

Let $m_1,m_2,\ldots,m_{g}$ be one-sided simple closed curves on $N_g$ as in Fig.~\ref{fig:one_sided_simple_closed_curves}.

\begin{figure}[ht]
    \begin{center}
        \includegraphics[scale=1.0,clip]{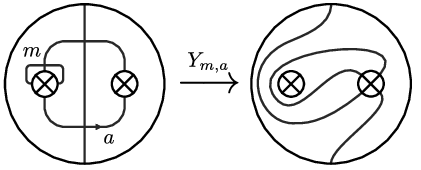}
    \end{center}
    \caption{Y-homeomorphism on $K$}
    \label{fig:y-homeo}
\end{figure}

%Y-homeoの定義
A Y-homeomorphism is defined as follow. For a one-sided simple closed curve $m$ and a two-sided oriented simple closed curve $a$ which intersects $m$ transversely in one point, the regular neighborhood $K$ of $m \cup a$ is homomeomorphic to the Klein bottle with one hole. Let $M$ be the regular neighborhood of $m$. Then the {\textit{Y-homeomorphism}} $Y_{m,a}$ is the isotopy class of the diffeomorphism obtained by pushing $M$ once along $a$ keeping the boundary of $K$ fixed (see Fig.~\ref{fig:y-homeo}).

Dehn twists and Y-homeomorphisms have the following properties.
%Dehn twistの性質%
\begin{lem}
    \label{lem:property_dehn_twist}
    For any element $f$ in $\CM _g$ and a two-sided simple closed curve $a$, we have

    \[ t_{f(a)}^{\epsilon} = f t_a f^{-1}, \]
    where if $f\mid_{N_a}$ is an orientation preserving diffeomorphism (resp.\ orientation reversing diffeomorphism), then $\epsilon = 1$ (resp. $\epsilon = -1$).
\end{lem}

%Y-homeomorphismの性質%
\begin{lem}
    \label{lemma:property_Y_homeomorphism}
    For a one-sided simple closed curve $m$, a two-sided simple closed curve $a$, and diffeomorphism $f$ on $N_g$, we have the following.
    \begin{align}
        Y_{m^{-1},a}     & = Y_{m,a},       \\
        Y_{m,a^{-1}}     & = Y_{m,a}^{-1},  \\
        f Y_{m,a} f^{-1} & = Y_{f(m),f(a)}.
    \end{align}
    %$(1)$ $Y_{m^{-1},a} = Y_{m,a}$. \\
    %$(2)$ $Y_{m,a^{-1}} = Y_{m,a}^{-1}$. \\
    %$(3)$ For any diffeomorphism $f$ of the surface $N_{g}$, we have $f Y_{m,a} f^{-1} = Y_{f(m),f(a)}$. \\
\end{lem}

%大森君の生成元を紹介
Omori \cite{o} reduced the number of Dehn twist generators for $\IT _g$ and showed:
\begin{thm}
    \label{thm:omori}
    For $g \geq 4$, $\IT _g$ is generated by $g+1$ Dehn twists $t_{a_1},t_{a_2},\ldots,t_{a_{g-1}},t_b,$ and $t_e$.
\end{thm}

%大森君の定理の注意を書く
\begin{rem}
    \label{rem:omori}
    {\rm{The minimum number of generators for $\IT _g$ by Dehn twists is at least $g$ for $g\geq4$ (see \cite{o}). Omori asked the following problem: {\it{Which of $g$ and $g+1$ is the minimum number of generators for $\IT _g$ by Dehn twists when $g\geq4$?}}}}
\end{rem}

\begin{figure}[ht]
    \begin{center}
        \includegraphics[scale=1.0,clip]{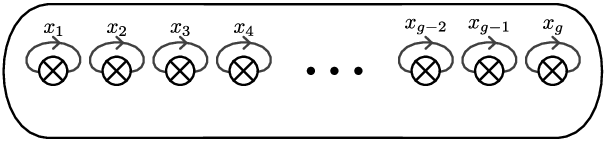}
    \end{center}
    \caption{A basis of ${\mathrm H_1}(N_g,\R)$.}
    \label{fig:homology_base}
\end{figure}

%以下ホモロジーについての話を書く ここはもっと自分の言葉にした方が良い
We set a basis %$x_1,x_2,\ldots, x_{g-1}$ 
of ${\mathrm H_1}(N_g,\R)$ as in Fig.~\ref{fig:homology_base}.

It is known that a homology group ${\mathrm H_1}(N_g,\R)$ has a basis in which every linear map $f_* : {\mathrm H_1}(N_g,\R) \rightarrow {\mathrm H_1}(N_g,\R)$, induced by a diffeomorphism $f:N_g \rightarrow N_g$, has a matrix with integral coefficients. Therefore we can define the homomorphism $D:\CM _g \rightarrow \Z _2$ as follows: $D(f) = det(f_*)$. We recall the properties of this homomorphism $D$ (see \cite{st1}).

%\begin{lem}
%	Let $c$ be a two-sided simple closed curve on $N_g$.
%	Then $D(t_c) = 1$.
%\end{lem}
%
%\begin{lem}
%	Let $a$ and $m$ be two-sided and one-sided simple closed curves on $N_g$, respectively. Supposed that $a$ intersects $m$ transversely in one point.
%	Then $D(Y_{m,a}) = -1$.
%\end{lem}
\begin{lem}
    \label{lem:prop_d}
    Let $f$ be an element of $\CM _g$.
    Then $D(f) = 1$ (resp.~$D(f) = -1$) if $f$ is (resp. is not) in $\IT _g$.
\end{lem}

\begin{rem}
    Let $y$ be any Y-homeomorphism. Then $D(y) = -1$.
\end{rem}

%------§4. Proof of Main Theorem 2 -----%
\section{Proof of Theorem~\ref{thm:main} }
\label{section:three}
%変形のところで
In this section, we prove Theorem~\ref{thm:main}. For $r\geq7$, We assume that $g = 2r$. %(In the case of $g=2r+1$, we can prove Theorem~\ref{thm:main} in same way).

\begin{figure}[ht]
    \begin{center}
        \includegraphics[scale=1.0,clip]{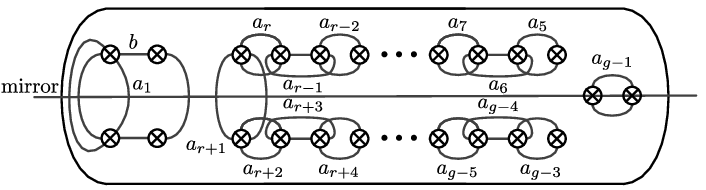}
    \end{center}
    \caption{Involution $\sigma : N_g \rightarrow N_g$.}
    \label{fig:involution_sigma_even}
\end{figure}

%\vspace{0.2cm}
%\begin{figure}[ht]
%\begin{center}
%\includegraphics[scale=1.0,clip]{./drawable/involution1}
%\end{center}
%\caption{Involution $\sigma : N_g \rightarrow N_g$.}
%\label{fig:involution_sigma}
%\end{figure}
%\clearpage
%
%\vspace{0.2cm}
%\begin{figure}[ht]
%\begin{center}
%\includegraphics[scale=1.0,clip]{./drawable/involution2}
%\end{center}
%\caption{Involution $\tau : N_g \rightarrow N_g$.}
%\label{fig:involution_tau}
%\end{figure}
%\clearpage
%
%\vspace{0.2cm}
%\begin{figure}[ht]
%\begin{center}
%\includegraphics[scale=1.0,clip]{./drawable/involution3}
%\end{center}
%\caption{Involution $\upsilon : N_g \rightarrow N_g$.}
%\label{fig:involution_upsilon}
%\end{figure}
%\clearpage
We deform $N_g$ in Fig.~\ref{fig:generator1} to $N_g$ in Fig.\ref{fig:involution_sigma_even} by diffeomorphism $\psi$ %which is isotopic to identify

so that $N_g$ in Fig~.\ref{fig:involution_sigma_even} is symmetrical with respect to the plane, denoted by mirror, illustrated in the middle of this surface. Let $\widetilde{\sigma}$ be a reflection in the mirror and let $\sigma$ be a composition $\psi^{-1} \widetilde{\sigma} \psi$. Then $\sigma$ is involution of $\CM _g$.
By an easy calculation, we found $D(\sigma) = -1$ if $r$ is odd and $D(\sigma) = 1$ if $r$ is even. When $g$ is odd, $\sigma$ is not element in $\IT _g$ by lemma \ref{lem:prop_d}. So, we consider a product $\sigma Y_{m_{g-1},a_{g-1}}$. This element $\sigma Y_{m_{g-1},a_{g-1}}$ is in $\IT _g$. Since $\sigma(a_{g-1}) = a_{g-1}$ and since $\sigma \mid _{N_{a_{g-1}}}$ is orientation reversing, we found that $\sigma Y_{m_{g-1},a_{g-1}}$ is an involution by lemma \ref{lemma:property_Y_homeomorphism}.%: $\sigma Y_{m_{g-1},a_{g-1}} \sigma Y_{m_{g-1},a_{g-1}} = Y_{\sigma(m_{g-1}),\sigma(a_{g-1})} Y_{m_{g-1},a_{g-1}} = Y_{m_{g-1}^{-1},a_{g-1}^{-1}} Y_{m_{g-1},a_{g-1}} = Y_{m_{g-1},a_{g-1}}^{-1} Y_{m_{g-1},a_{g-1}} = id$. 
We rewrite $\sigma Y_{m_{g-1},a_{g-1}}$ as $\sigma$.
Therefore, for $r\geq7$, $\sigma$ is involution in $\IT _g$.

\begin{figure}[ht]
    \begin{center}
        \includegraphics[scale=1.0,clip]{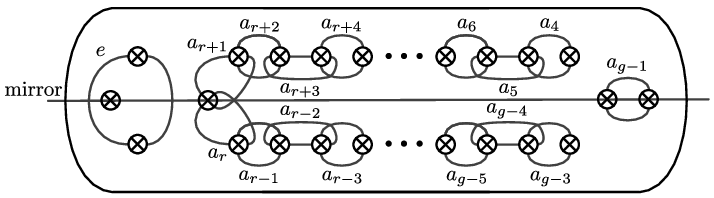}
    \end{center}
    \caption{Involution $\tau : N_g \rightarrow N_g$.}
    \label{fig:involution_tau_even}
\end{figure}

In the same way, we define $\tau$ as reflection in the mirror in Fig.~\ref{fig:involution_tau_even}.

We found $D(\tau) = 1$ if $r$ is odd and $D(\tau) = -1$ if $r$ is even. When $g$ is even, $\tau$ is not element in $\IT _g$ by lemma \ref{lem:prop_d}. So, we consider a product $\tau Y_{m_{g-1},a_{g-1}}$. This element $\tau Y_{m_{g-1},a_{g-1}}$ is in $\IT _g$. By lemma \ref{lemma:property_Y_homeomorphism}, $\tau Y_{m_{g-1},a_{g-1}}$ is an involution. We rewrite $\tau Y_{m_{g-1},a_{g-1}}$ as $\tau$.
Therefore, for $r\geq7$, $\tau$ is involution in $\IT _g$.

\begin{figure}[ht]
    \begin{center}
        \includegraphics[scale=1.0,clip]{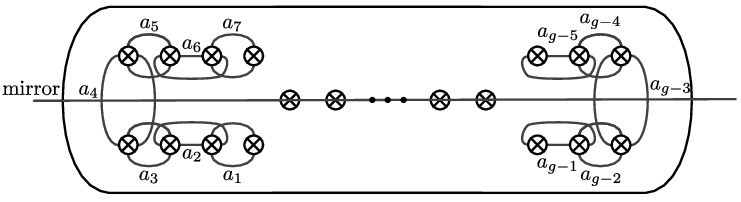}
    \end{center}
    \caption{Involution $\upsilon : N_g \rightarrow N_g$.}
    \label{fig:involution_upsilon_even}
\end{figure}

We construct third involution. We define $\upsilon$ as reflection in the mirror in Fig.~\ref{fig:involution_upsilon_even} and we found $D(\upsilon) = 1$. Hence, the element $\upsilon$ is involution in $\IT _g$.

Let $\rho _1$, $\rho _2$, and $\rho _3$ be $\sigma t_{a_1}$, $\sigma t_{b}$, and $\tau t_{e}$, respectively. By lemma~\ref{lem:property_dehn_twist}, $\rho _i$ is involution in $\IT _g$ for $i=1,2,3$. Let $G$ be a subgroup of $\IT _g$ generated by $\sigma$, $\tau$, $\upsilon$, $\rho _1$, $\rho _2$, and $\rho _3$. Since $t_{a_1}$ is a product $\sigma \rho _1$, $t_{a_1}$ is a product of two involutions. In the same way, $t_b$ and $t_e$ are products of two involutions. We show that $t_{a_2},t_{a_3},\ldots,t_{a_{g-1}}$ are products of involutions.

\begin{figure}[ht]
    \begin{center}
        \includegraphics[scale=1.0,clip]{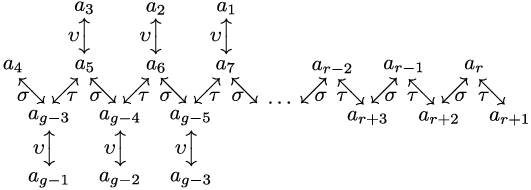}
    \end{center}
    \caption{}
    \label{fig:action}
\end{figure}

Let $a$ and $b$ be two-sided simple closed curves on $N_g$. For $f \in G$, the symbol $a \stackrel{f}{\longleftrightarrow} b$ means that $f(a) = b$ or $f^{-1}(a) = b$.
%Since $t_{a_1}$ is mapped to $t_{a_2},t_{a_3},\ldots,t_{a_{g-1}}$ (see Fig.~\ref{fig:action}), $t_{a_i}$ is a product of involutions by lemma~\ref{lem:property_dehn_twist} for $i=1,2,\ldots,g-1$.
Since $a_1$ is mapped to $a_2,a_3,\ldots,a_{g-1}$ by $\sigma$, $\tau$, and $\upsilon$ (see Fig.~\ref{fig:action}), $t_{a_i}$ is a product of involutions by lemma~\ref{lem:property_dehn_twist} for $i=1,2,\ldots,g-1$.
Since all Omori's generators for $\IT _g$ is in $G$, we prove that $G$ is equal to $\IT _g$.

\begin{figure}[ht]
    \begin{center}
        \includegraphics[scale=1.0,clip]{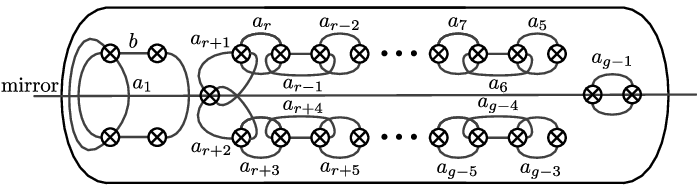}
    \end{center}
    \caption{Involution $\sigma : N_g \rightarrow N_g$.}
    \label{fig:involution_sigma_odd}
\end{figure}

\begin{figure}[ht]
    \begin{center}
        \includegraphics[scale=1.0,clip]{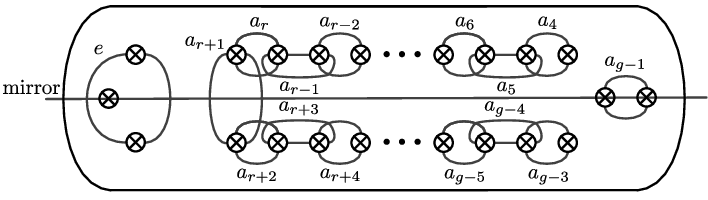}
    \end{center}
    \caption{Involution $\tau : N_g \rightarrow N_g$.}
    \label{fig:involution_tau_odd}
\end{figure}

\begin{figure}[ht]
    \begin{center}
        \includegraphics[scale=1.0,clip]{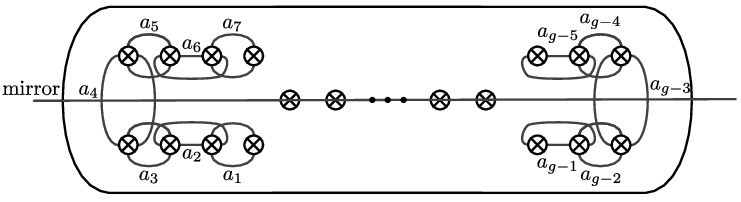}
    \end{center}
    \caption{Involution $\upsilon' : N_g \rightarrow N_g$.}
    \label{fig:involution_upsilon_odd}
\end{figure}

In the case of $g = 2r+1 (r\geq7)$, we can prove Theorem~\ref{thm:main} in same way. Involutions $\sigma$, $\tau$, and $\upsilon'$ are reflection in the mirror in Fig.~\ref{fig:involution_sigma_odd} , Fig.~\ref{fig:involution_tau_odd} , and Fig.~\ref{fig:involution_upsilon_odd} , respectively.

If $\sigma$ (resp. $\tau$) is not in $\IT _g$, we rewrite $\sigma Y_{m_{g-1},a_{g-1} }$ (resp.~ $\tau Y_{m_{g-1},a_{g-1} }$) as $\sigma$ (resp. $\tau$). Since the element $\upsilon'$ is not in $\IT _g$, we define $\upsilon$ as a product $\upsilon' Y_{m_9,a_9}$. This element $\upsilon$ is in $\IT _g$ and this element is an involution by lemma \ref{lemma:property_Y_homeomorphism}. Involution $\rho_1$, $\rho_2$, and $\rho_3$ are the same as the case of $g=2r$. %Since $\t_{a_1}$ is a prodcut of two involutions and since $a_1$ is mapped to $a_2,a_3,\ldots,a_{g-1}$ by $\sigma$, $\tau$, and $\upsilon$, $t_{a_i}$ is a product of involutions for $i=1,2,\ldots,g-1$.
In the same way as the case of $g=2r$, all Omori's generators for $\IT_g$ is in $G$.

So we proved Theorem~\ref{thm:main2}:
\begin{thm}
    \label{thm:main2}
    For $g\geq16$ or $g=14$, $\IT_g$ is generated by involutions $\sigma$, $\tau$, $\upsilon$, $\rho_1$, $\rho_2$, and $\rho_3$.
\end{thm}

%------§4. Remark -----%
\section{Concluding Remarks}
\label{section:four}
\begin{figure}[ht]
    \begin{center}
        \includegraphics[scale=1.0,clip]{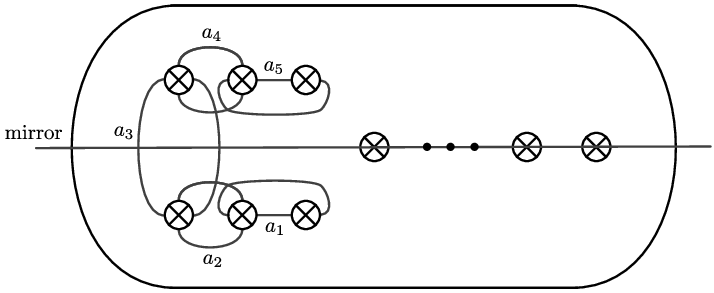}
    \end{center}
    \caption{Involution $\eta : N_g \rightarrow N_g$.}
    \label{fig:involution_eta}
\end{figure}

\begin{figure}[ht]
    \begin{center}
        \includegraphics[scale=1.0,clip]{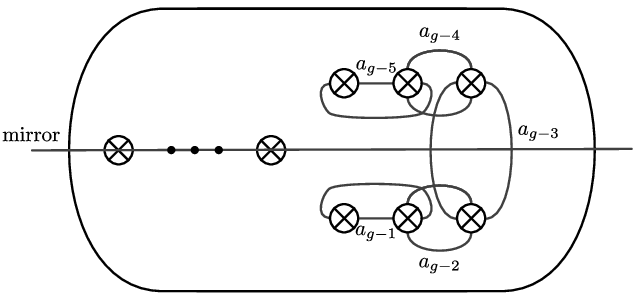}
    \end{center}
    \caption{Involution $\theta : N_g \rightarrow N_g$.}
    \label{fig:involution_theta}
\end{figure}

We found the construction of $\upsilon$ needs $g\geq 16$ or $g=14$. By constructing of three other involutions instead of $\upsilon$, $\IT _g$ can be generated by eight involutions in case of $g \geq 8$.
In this case, $\IT _g$ is generated by $\sigma,\tau,\rho_1,\rho_2,\rho_3$ in Section $3$ and additionally $\eta$, $\theta$, and $\rho _4$. Let $\eta'$ and $\theta'$ are involutions of $\IT _g$ which are reflections in the mirrors in Fig.~\ref{fig:involution_eta} and Fig.~\ref{fig:involution_theta}, respectively. We define $\eta = \eta' Y_{m_7,a_7}$ and $\theta = \theta' Y_{m_1,a_1}$. A element $\rho _4$ is a product of $\eta t_{a_3}$. Then, $\eta$, $\theta$, and $\rho _4$ are involutions in $\IT _g$ by lemma \ref{lemma:property_Y_homeomorphism} and \ref{lem:prop_d}. In the same way as Section 3, We can proved the following:

\begin{thm}
    \label{thm:main_remark}
    For $g \geq 8$, $\IT _g$ is generated by eight involutions.
\end{thm}

\vspace{0.2cm}

\vspace{0.2cm}

In the case of $\CM _g$, Szepietowski show that $\CM _g$ is generated by $4$ involutions for $g \geq 4$.

%\begin{thm}
%	$\CM _g$ is generated by $4$ involutions for $g \geq 4$.
%\end{thm}
We can consider following problem:
\begin{prob}
    For $g\geq4$, can the twist subgroup $\IT _g$ be generated by $4$ involutions?
\end{prob}

\section*{Acknowledgments}
The author would like to thank Susumu Hirose and Naoyuki Monden for many advices. He would also like to thank Genki Omori for many helpful suggestions and comments of the twist subgroup.

%---- Reference ----%

% TODOを書き込む

\end{document}